\documentclass[reqno,12pt]{amsart}
\usepackage{amssymb}
\usepackage{amsmath}
\usepackage{amsfonts}

\newtheorem{theorem}{Theorem}[section]
\theoremstyle{plain}

\newtheorem{corollary}[theorem]{Corollary}

\newtheorem{definition}[theorem]{Definition}

\newtheorem{lemma}[theorem]{Lemma}

\newtheorem{remark}[theorem]{Remarque}

\numberwithin{equation}{section}
\input{tcilatex}

\begin{document}
\title{ Caracterization of Fr\'{e}chet differentiability norm for dual
complex Banach spaces}
\author{Mohammad Daher}
\email{daher.mohammad@ymail.com}

\begin{abstract}
Let $X$ be a complex Banach space; in this work we characterize the property
of Fr\'{e}chet differentiability for the dual space of $X.$ In the
following, we show that if the dual space of $X$ is G\^{a}teaux
differentiable, then the dual space of $L^{p}(X)$ has the same property for
all $1<p<+\infty .$
\end{abstract}

\maketitle

0\bigskip AMS Classification : 46B70

Mots cl\'{e}s : Interpolation, smooth\ \ \ \ \ \ \ \ \ \ \ \ \ \ \ \ \ \ \ \
\ \ \ \ \ \ \ \ \ \ \ \ \ \ \ \ \ \ \ \ \ \ \ \ \ \ \ \ \ \ \ \ \ \ \ \ \ \
\ \ \ \ \ \ \ \ \ \ \ \ \ \ \ \ \ \ \ \ \ \ \ \ \ \ \ \ \ \ \ \ \ \ \ \ \ \
\ \ \ \ \ \ \ \ \ \ \ \ \ \ \ \ \ \ \ \ \ \ \ \ \ \ \ \ \ \ \ \ \ \ \ \ \ \
\ \ \ \ \ \ \ \ \ \ \ \ \ \ \ \ \ \ \ \ \ \ \ \ \ \ \ \ \ \ \ \ \ \ \ \ \ \
\ \ \ \ \ \ \ \ \ \ \ \ \ \ \ \ \ \ \ \ \ \ \ \ \ \ \ \ \ \ \ \ \ \ \ \ \ \
\ \ \ \ \ \ \ \ \ \ \ \ \ \ \ \ \ \ \ \ \ \ \ \ \ \ \ \ \ \ \ \ \ \ \ \ \ \
\ \ \ \ \ \ \ \ \ \ \ \ \ \ \ \ \ \ \ \ \ \ \ \ \ \ \ \ \ \ \ \ \ \ \ \ \ \
\ \ \ \ \ \ \ \ \ \ \ \ \ \ \ \ \ \ \ \ \ \ \ \ \ \ \ \ \ \ \ \ \ \ \ \ \ \
\ \ \ \ \ \ \ \ \ \ \ \ \ \ \ \ \ \ \ \ \ \ \ \ \ \ \ \ \ \ \ \ \ \ \ \ \ \
\ \ \ \ \ \ \ \ \ \ \ \ \ \ \ \ \ \ \ \ \ \ \ \ \ \ \ \ \ \ \ \ \ \ \ \ \ \
\ \ \ \ \ \ \ \ \ \ \ \ \ \ \ \ \ \ \ \ \ \ \ \ \ \ \ \ \ \ \ \ \ \ \ \ \ \
\ \ \ \ \ \ \ \ \ \ \ \ \ \ \ \ \ \ \ \ \ \ \ \ \ \ \ \ \ \ \ \ \ \ \ \ \ \
\ \ \ \ \ \ \ \ \ \ \ \ \ \ \ \ \ \ \ \ \ \ \ \ \ \ \ \ \ \ \ \ \ \ \ \ \ \
\ \ \ \ \ \ \ \ \ \ \ \ \ \ \ \ \ \ \ \ \ \ \ \ \ \ \ \ \ \ \ \ \ \ \ \ \ \
\ \ \ \ \ \ \ \ \ \ \ \ \ \ \ \ \ \ \ \ \ \ \ \ \ \ \ \ \ \ \ \ \ \ \ \ \ \
\ \ \ \ \ \ \ \ \ \ \ \ \ \ \ \ \ \ \ \ \ \ \ \ \ \ \ \ \ \ \ \ \ \ \ \ \ \
\ \ \ \ \ \ \ \ \ \ \ \ \ \ \ \ \ \ \ \ \ \ \ \ \ \ \ \ \ \ \ \ \ \ \ \ \ \
\ \ \ \ \ \ \ \ \ \ \ \ \ \ \ \ \ \ \ \ \ \ \ \ \ \ \ \ \ \ \ \ \ \ \ \ \ \
\ \ \ \ \ \ \ \ \ \ \ \ \ \ \ \ \ \ \ \ \ \ \ \ \ \ \ \ \ \ \ \ \ \ \ \ \ \
\ \ \ \ \ \ \ \ \ \ \ \ \ \ \ \ \ \ \ \ \ \ \ \ \ \ \ \ \ \ \ \ \ \ \ \ \ \
\ \ \ \ \ \ \ \ \ \ \ \ \ \ \ \ \ \ \ \ \ \ \ \ \ \ \ \ \ \ \ \ \ \ \ \ \ \
\ \ \ \ \ \ \ \ \ \ \ \ \ \ \ \ \ \ \ \ \ \ \ \ \ \ \ \ \ \ \ \ \ \ \ \ \ \
\ \ \ \ \ \ \ \ \ \ \ \ \ \ \ \ \ \ \ \ \ \ \ \ \ \ \ \ \ \ \ \ \ \ \ \ \ \
\ \ \ \ \ \ \ \ \ \ \ \ \ \ \ \ \ \ \ \ \ \ \ \ \ \ \ \ \ \ \ \ \ \ \ \ \ \
\ \ \ \ \ \ \ \ \ \ \ \ \ \ \ \ \ \ \ \ \ \ \ \ \ \ \ \ \ \ \ \ \ \ \ \ \ \
\ \ \ \ \ \ \ \ \ \ \ \ \ \ \ \ \ \ \ \ \ \ \ \ \ \ \ \ \ \ \ \ \ \ \ \ \ \
\ \ \ \ \ \ \ \ \ \ \ \ \ \ \ \ \ \ \ \ \ \ \ \ \ \ \ \ \ \ \ \ \ \ \ \ \ \
\ \ \ \ \ \ \ \ \ \ \ \ \ \ \ \ \ \ \ \ \ \ \ \ \ \ \ \ \ \ \ \ \ \ \ \ \ \
\ \ \ \ \ \ \ \ \ \ \ \ \ \ \ \ \ \ \ \ \ \ \ \ \ \ \ \ \ \ \ \ \ \ \ \ \ \
\ \ \ \ \ \ \ \ \ \ \ \ \ \ \ \ \ \ \ \ \ \ \ \ \ \ \ \ \ \ \ \ \ \ \ \ \ \
\ \ \ \ \ \ \ \ \ \ \ \ \ \ \ \ \ \ \ \ \ \ \ \ \ \ \ \ \ \ \ \ \ \ \ \ \ \
\ \ \ \ \ \ \ \ \ \ \ \ \ \ \ \ \ \ \ \ \ \ \ \ \ \ \ \ \ \ \ \ \ \ \ \ \ \
\ \ \ \ \ \ \ \ \ \ \ \ \ \ \ \ \ \ \ \ \ \ \ \ \ \ \ \ \ \ \ \ \ \ \ \ \ \
\ \ \ \ \ \ \ \ \ \ \ \ \ \ \ \ \ \ \ \ \ \ \ \ \ \ \ \ \ \ \ \ \ \ \ \ \ \
\ \ \ \ \ \ \ \ \ \ \ \ \ \ \ \ \ \ \ \ \ \ \ \ \ \ \ \ \ \ \ \ \ \ \ \ \ \
\ \ \ \ \ \ \ \ \ \ \ \ \ \ \ \ \ \ \ \ \ \ \ \ \ \ \ \ \ \ \ \ \ \ \ \ \ \
\ \ \ \ \ \ \ \ \ \ \ \ \ \ \ \ \ \ \ \ \ \ \ \ \ \ \ \ \ \ \ \ \ \ \ \ \ \
\ \ \ \ \ \ \ \ \ \ \ \ \ \ \ \ \ \ \ \ \ \ \ \ \ \ \ \ \ \ \ \ \ \ \ \ \ \
\ \ \ \ \ \ \ \ \ \ \ \ \ \ \ \ \ \ \ \ \ \ \ \ \ \ \ \ \ \ \ \ \ \ \ \ \ \
\ \ \ \ \ \ \ \ \ \ \ \ \ \ \ \ \ \ \ \ \ \ \ \ \ \ \ \ \ \ \ \ \ \ \ \ \ \
\ \ \ \ \ \ \ \ \ \ \ \ \ \ \ \ \ \ \ \ \ \ \ \ \ \ \ \ \ \ \ \ \ \ 

\section{\textsc{Introduction}}

It is known the criteria of the property of Fr\'{e}chet differentiability
for the real Banach spaces (see \cite[Chap.I]{D-G-Z}). In this work we give
a similar criteria of the property of Fr\'{e}chet differentiability for the
dual complex Banach spaces. In the rest of this work we show that if $%
X^{\ast }$ is G\^{a}teaux differentiable, then $L^{p}(\Omega ,X)^{\ast }$ is
G\^{a}teaux differentiable for all $1<p<+\infty $.

In \cite{Le-Su} one shows that if the Banach space $X$ is Fr\'{e}chet
differentiable, then $L^{p}(\Omega ,X)$ has the same propertyfor all $%
1<p<+\infty $ . In the reste of this work, we show that if $X^{\ast }$ is G%
\^{a}teaux differentiable, then $L^{p}(\Omega ,X)^{\ast }$ is G\^{a}teaux
differentiable for all $1<p<+\infty $, where $(\Omega ,\mu )$ is a
probability space.

Let $X$ be a Banach space. We denote by $B_{X}$ the closed unit ball of a
Banach space $X$ and by $S_{X}$ its unit sphere. On the other hand we denote
by $\left\langle x,x^{\ast }\right\rangle $ the pairing between an element $%
x $ of $X$ and an element $x^{\ast }$ of its dual $X^{\ast }$.

\begin{definition}
\label{tz}A Banach space $(Y,.\bigl\Vert.\bigr\Vert)$ is said to be Fr\'{e}%
chet differentiable, if it is Fr\'{e}chet differentiable at every point $%
a\in Y$, $a\neq 0$, \textit{i.e.}, lim$_{\bigl\Vert h\bigr\Vert\rightarrow 0}%
\frac{\bigl\Vert a+h\bigr\Vert+\bigl\Vert a-h\bigr\Vert-2\bigl\Vert a%
\bigr\Vert}{\bigl\Vert h\bigr\Vert}=0.$
\end{definition}

\begin{definition}
A Banach space $(Y,\bigl\Vert.\bigr\Vert)$ is said to be G\^{a}teaux
differentiable, if it is G\^{a}teaux differentiable at every point $a\in Y$, 
$a\neq 0$, \textit{i.e.}, if for every $h\in Y$, $\lim_{\mathbb{R}\ni
t\rightarrow 0}\frac{\bigl\Vert a+th\bigr\Vert-\bigl\Vert a\bigr\Vert}{t}$
exists.
\end{definition}

\begin{definition}
\label{Sm copy(1)} A Banach space $(Y,\bigl\Vert.\bigr\Vert)$ is said to be
smooth, if it is smooth at every point $a\in Y$, $a\neq 0$, \textit{i.e.},
there exists a unique $z^{\ast }\in S_{Y^{\ast }}$ such that $\bigl\Vert a%
\bigr\Vert=\left\langle a,z^{\ast }\right\rangle $.
\end{definition}

\begin{definition}
\bigskip \label{UL} A Banach space $(Y,\bigl\Vert.\bigr\Vert)$ is said to be
uniformly smooth, if for every $x\in S_{X},$ there is a unique point $%
x^{\ast }\in S_{X^{\ast }}$ satisfying for every $\varepsilon >0,$ there
exists $\delta >0$ such that if $g\in B_{X^{\ast }}$ with $\left\langle
x,g\right\rangle >1-\delta ,$ then $\bigl\Vert x^{\ast }-g\bigr\Vert%
<\varepsilon $ $.$

Let $X$ be a real Banach space; by the \v{S}mulyan's result \cite[%
Chap.I,Th.1.4,Cor.1.5]{D-G-Z} the dual space $X^{\ast }$ is Fr\'{e}chet
differentiable if and only if $X^{\ast }$ satisfying one of the following
conditions:
\end{definition}

(1) $X^{\ast }$ is uniformly smooth.

(2) For all sequences $(x_{n})_{n\geq 0},(y_{n})_{n\geq 0}$ in $S_{X}$ and
all $a^{\ast }\in S_{X}$ such that $\left\langle x_{n},a^{\ast
}\right\rangle \rightarrow _{n\rightarrow +\infty }1$ and $\left\langle
y_{n},a^{\ast }\right\rangle \rightarrow _{n\rightarrow +\infty }1$, then $%
\bigl\Vert x_{n}-y_{n}\bigr\Vert_{X}\rightarrow _{n\rightarrow +\infty }0.$

Recall that by \cite[Lemma 5.2]{Da1} if $X$ is a complex Banach, then $%
X^{\ast }$ is Fr\'{e}chet differentiable if and only if $X^{\ast }$
satisfies the condition (2). The proof for the complex case is exactly
similar to the proof in \cite[Chap.I,Th.1.4]{D-G-Z} with a small
modification.

\begin{theorem}
\label{K}Let $X$ be a complex Banach space. Then $X^{\ast }$ is uniformly
smooth if and only if $X^{\ast }$ satisfies the condition 2).
\end{theorem}

Proof: Assume that $X^{\ast }$ is uniformly smooth. Let $(f_{n})_{n\geq
0},(g_{n})_{n\geq 0}$ be two sequences in $S_{X}$ and let $x^{\ast }\in
S_{X^{\ast }}$ such that $\left\langle f_{n},x^{\ast }\right\rangle
\rightarrow _{n\rightarrow +\infty }1$ and $\left\langle g_{n},x^{\ast
}\right\rangle \rightarrow _{n\rightarrow +\infty }1.$ Since $X^{\ast }$ is
uniformly smooth, there is $x^{\ast \ast }$ a unique point in $S_{X^{\ast
\ast }}$ satisfying for every $\varepsilon >0,$ there exists $\delta >0$
such that if $y^{\ast \ast }\in B_{X^{\ast \ast }}$ with $\left\langle
x^{\ast },y^{\ast \ast }\right\rangle >1-\delta ,$ then $\bigl\Vert x^{\ast
\ast }-y^{\ast \ast }\bigr\Vert<\frac{\varepsilon }{2}.$ Let $\varepsilon >0$%
, put $u_{n}=\frac{f_{n}}{\left\langle f_{n},x^{\ast }\right\rangle }%
\left\vert \left\langle f_{n},x^{\ast }\right\rangle \right\vert $ and $%
v_{n}=\frac{g_{n}}{\left\langle g_{n},x^{\ast }\right\rangle }\left\vert
\left\langle g_{n},x^{\ast }\right\rangle \right\vert .$ As $\lambda _{n}=%
\frac{\left\vert \left\langle f_{n},x^{\ast }\right\rangle \right\vert }{%
\left\langle f_{n},x^{\ast }\right\rangle }\rightarrow _{n\rightarrow
+\infty }1$ and $\mu _{n}=\frac{\left\vert \left\langle g_{n},x^{\ast
}\right\rangle \right\vert }{\left\langle g_{n},x^{\ast }\right\rangle }%
\rightarrow _{n\rightarrow +\infty }1,$ $\left\langle u_{n},x^{\ast
}\right\rangle \rightarrow _{n\rightarrow +\infty }1$ and $\left\langle
v_{n},x^{\ast }\right\rangle \rightarrow _{n\rightarrow +\infty }1.$ But, $%
\left\langle u_{n},x^{\ast }\right\rangle ,$ $\left\langle v_{n},x^{\ast
}\right\rangle \in \mathbb{R}^{+},$ hence there is $n_{0}\in \mathbb{N}$
such that $\left\langle u_{n},x^{\ast }\right\rangle >1-\delta $ and $%
\left\langle v_{n},x^{\ast }\right\rangle >1-\delta ,$ for all $n\geq n_{0};$
it follows that $\bigl\Vert x^{\ast \ast }-u_{n}\bigr\Vert<\frac{\varepsilon 
}{2}$ and $\bigl\Vert x^{\ast \ast }-v_{n}\bigr\Vert<\frac{\varepsilon }{2},$
this implies that $\bigl\Vert u_{n}-v_{n}\bigr\Vert<\varepsilon .$ Thus $%
\bigl\Vert u_{n}-v_{n}\bigr\Vert\rightarrow _{n\rightarrow +\infty }0.$ On
the other hand $f_{n}=\frac{u_{n}}{\lambda _{n}}$, $g_{n}=\frac{v_{n}}{\mu
_{n}},$ $\lambda _{n}\rightarrow _{n\rightarrow +\infty }1$ and $\mu
_{n}\rightarrow _{n\rightarrow +\infty }1,$ we deduce that $\bigl\Vert %
f_{n}-g_{n}\bigr\Vert\rightarrow _{n\rightarrow +\infty }0.$\qed

Conversely, assume that $X^{\ast }$ is not uniformly smooth. So there exist $%
a^{\ast }\in X^{\ast },$ $a^{\ast \ast }\in X^{\ast \ast },$ $\varepsilon >0$
and a sequence $(x_{n}^{\ast \ast })_{n\geq 0}$ in $B_{X^{\ast \ast }}$ such
that $\left\langle a^{\ast },a^{\ast \ast }\right\rangle =1,$ $\bigl\Vert %
x_{n}^{\ast \ast }-a^{\ast \ast }\bigr\Vert>\varepsilon $ and \ 
\begin{equation}
\left\langle a^{\ast },x_{n}^{\ast \ast }\right\rangle >1-\frac{1}{n+1},
\label{DF}
\end{equation}

for all $n\geq 0.$ Now observe that for every $n\geq 0,$ there is $%
u_{n}^{\ast }$ in the unit ball of $X^{\ast }$ such that%
\begin{equation}
\left\langle u_{n},x_{n}^{\ast \ast }-a^{\ast \ast }\right\rangle
>\varepsilon .  \label{LL}
\end{equation}%
For every $n\geq 0$ denote 
\begin{eqnarray*}
D_{n} &=&\left\{ x^{\ast \ast }\in B_{X^{\ast \ast }};\text{ }\left\vert
\left\langle a^{\ast },x^{\ast \ast }\right\rangle -1\right\vert <\frac{1}{%
n+1}\right\} \cap \\
&&\left\{ x^{\ast \ast }\in B_{X^{\ast \ast }};\text{ }\left\vert
\left\langle u_{n}^{\ast },x^{\ast \ast }\right\rangle -\left\langle
u_{n}^{\ast },a^{\ast \ast }\right\rangle \right\vert <\frac{1}{n+1}\right\}
\end{eqnarray*}%
and 
\begin{eqnarray*}
E_{n} &=&\left\{ x^{\ast \ast }\in B_{X^{\ast \ast }};\text{ }\left\vert
\left\langle a^{\ast },x_{{}}^{\ast \ast }\right\rangle -\left\langle
a^{\ast },x_{n}^{\ast \ast }\right\rangle \right\vert <\frac{1}{n+1}\right\}
\cap \\
&&\left\{ x^{\ast \ast }\in B_{X^{\ast \ast }};\text{ }\left\vert
\left\langle u_{n}^{\ast },x_{{}}^{\ast \ast }\right\rangle -\left\langle
u_{n}^{\ast },x_{n}^{\ast \ast }\right\rangle \right\vert <\frac{1}{n+1}%
\right\} .
\end{eqnarray*}

Since $D_{n}$ is a $w^{\ast }-$neighborhood of $a^{\ast \ast }$ and $E_{n}$
is a $w^{\ast }-$neighborhood of $x_{n}^{\ast \ast },$ there exist $x_{n}\in
D_{n}\cap X$ and $y_{n}\in E_{n}\cap X.$ Thus%
\begin{equation}
\left\langle x_{n},a^{\ast }\right\rangle \rightarrow _{n\rightarrow +\infty
}1,  \label{X}
\end{equation}%
\begin{equation}
\left\langle x_{n},u_{n}^{\ast }\right\rangle -\left\langle u_{n}^{\ast
},a^{\ast \ast }\right\rangle \rightarrow _{n\rightarrow +\infty }0,
\label{B}
\end{equation}

\begin{equation}
\left\langle y_{n},a^{\ast }\right\rangle -\left\langle a^{\ast
},x_{n}^{\ast \ast }\right\rangle \rightarrow _{n\rightarrow +\infty }0,
\label{KU}
\end{equation}

and%
\begin{equation}
\left\langle y_{n},u_{n}^{\ast }\right\rangle -\left\langle \left\langle
u_{n}^{\ast },x_{n}^{\ast \ast }\right\rangle \right\rangle \rightarrow
_{n\rightarrow +\infty }0.  \label{H}
\end{equation}

On the other hand, by (\ref{DF}), $\left\langle a^{\ast },x_{n}^{\ast \ast
}\right\rangle \rightarrow _{n\rightarrow +\infty }1;$ this implies by (\ref%
{KU}) that%
\begin{equation}
\left\langle y_{n},a^{\ast }\right\rangle \rightarrow _{n\rightarrow +\infty
}1.  \label{O}
\end{equation}

Using (\ref{B}) and (\ref{H}) we find $n_{0}\in \mathbb{N}$ such that $%
\left\vert \left\langle x_{n},u_{n}^{\ast }\right\rangle -\left\langle
u_{n}^{\ast },a^{\ast \ast }\right\rangle \right\vert <\frac{\varepsilon }{4}
$ and $\left\vert \left\langle u_{n}^{\ast },x_{n}^{\ast \ast }\right\rangle
-\left\langle y_{n},u_{n}^{\ast }\right\rangle \right\vert <\frac{%
\varepsilon }{4},$ for all $n\geq n_{0}.$ It results by(\ref{LL}) that 
\begin{eqnarray*}
\bigl\Vert x_{n}-y_{n}\bigr\Vert &\geq &\left\vert \left\langle
x_{n}-y_{n},u_{n}^{\ast }\right\rangle \right\vert \\
&\geq &\left\vert \left\langle u_{n}^{\ast },a^{\ast \ast }\right\rangle
-\left\langle u_{n}^{\ast },x_{n}^{\ast \ast }\right\rangle \right\vert
-\left\vert \left\langle x_{n},u_{n}^{\ast }\right\rangle -\left\langle
u_{n}^{\ast },a^{\ast \ast }\right\rangle \right\vert -\left\vert
\left\langle u_{n}^{\ast },x_{n}^{\ast \ast }\right\rangle -\left\langle
y_{n},u_{n}^{\ast }\right\rangle \right\vert \geq \frac{\varepsilon }{2},
\end{eqnarray*}

for all $n\geq n_{0}.$ Since $\left\langle x_{n},a^{\ast }\right\rangle
\rightarrow _{n\rightarrow +\infty }1$ and $\left\langle y_{n},a^{\ast
}\right\rangle \rightarrow _{n\rightarrow +\infty }1,$ by (\ref{X}) and (\ref%
{O})$,$ $X^{\ast }$ does not satisfie the property 2. \qed

By theorem \ref{K} and \cite[lemma 5.2]{Da} one has the following result:

\begin{corollary}
\label{wq}Let $X$ be a complex Banach space. The following assertions are
equivalent:
\end{corollary}

1) $X^{\ast }$ is Fr\'{e}chet differentiable.

2) For all sequences $(x_{n})_{n\geq 0},(y_{n})_{n\geq 0}$ in $S_{X}$ and
all $a^{\ast }\in S_{X^{\ast }}$ such that $\left\langle x_{n},a^{\ast
}\right\rangle \rightarrow _{n\rightarrow +\infty }1$ and $\left\langle
y_{n},a^{\ast }\right\rangle \rightarrow _{n\rightarrow +\infty }1$, we have 
$\bigl\Vert x_{n}-y_{n}\bigr\Vert_{X}\rightarrow _{n\rightarrow +\infty }0.$

3) $X^{\ast }$ is uniformly smooth.

\begin{remark}
\label{GT}Let $X$ be a complex Banach space. Then $X$ is uniformly smooth if
and only if for all sequences $(f_{n}^{\ast })_{n\geq 0},(g_{n}^{\ast
})_{n\geq 0}$ in $S_{X^{\ast }}$ and all $x\in S_{X}$ such that $%
\left\langle x,f_{n}^{\ast }\right\rangle \rightarrow _{n\rightarrow +\infty
}1$ and $\left\langle x,g_{n}^{\ast }\right\rangle \rightarrow
_{n\rightarrow +\infty }1,$ then $\bigl\Vert f_{n}^{\ast }-g_{n}^{\ast }%
\bigr\Vert\rightarrow _{n\rightarrow +\infty }0$ in the weak topology$.$ The
proof of this remark is similar to that one of \cite[Chap.I,Th.1.4]{D-G-Z}.
\end{remark}

Let $(\Omega ,\mu )$ be a probability space, let $1<p<+\infty $ and let $X$
be a seperable Banach space$;$ denote

\begin{equation*}
VB^{p}(\Omega ,X^{\ast })=\left\{ f:\Omega \rightarrow X^{\ast }\text{ }%
w^{\ast }-measurable\text{ and }\mathop{\displaystyle \int}\limits_{\Omega }%
\bigl\Vert f(\omega )\bigr\Vert_{X^{\ast }}^{p}d\mu (\omega )<+\infty
\right\} .
\end{equation*}

\begin{theorem}
\label{ts}Let $(\Omega ,\mu )$ be a probabilty space and let $X$ be a
Banach. Assume that $X^{\ast }$ is G\^{a}teaux differentiable. Then $%
L^{p}(\Omega ,X)^{\ast }$ is G\^{a}teaux differentiable, for every $%
1<p<+\infty $.
\end{theorem}

Before giving the proof of Theorem \ref{ts}, We need the following lemmas .

\begin{lemma}
\label{vvv}\cite[Lemma 4]{Da2}A Banach space $(X,\Vert \mkern1.0mu.\mkern%
1.0mu\Vert )$ is G\^{a}teaux differentiable if for some $p\in \lbrack
1,+\infty )$ (resp. if for \emph{every} $p\in \lbrack 1,+\infty )$) we have
that:

for all $x\in S_{X}$ and $h\in X$, 
\begin{equation*}
\lim_{\mathbb{R}\ni t\rightarrow 0}\frac{\Vert x+th\Vert ^{p}+\Vert
x-th\Vert ^{p}-2\Vert x\Vert ^{p}}{t}=0.
\end{equation*}
\end{lemma}

\begin{lemma}
\label{UN}\cite[Lemma 5]{Da2}Let $p$ satisfy $1\leqslant p<+\infty $. For
every normed space~$X$, for all $x,h\in X$ we have that 
\begin{equation*}
\max_{t\in (0,1]}\frac{\Vert x+th\Vert ^{p}+\Vert x-th\Vert ^{p}-2\Vert
x\Vert ^{p}}{t}\leqslant 2^{p}(\Vert x\Vert ^{p}+\Vert h\Vert ^{p}).
\end{equation*}
\end{lemma}

\noindent

\bigskip Proof of theorem \ref{ts}: Let $Y$ be a closed separable subspace
of $L^{p}(\Omega ,X);$ there is a closed separable subspace $X_{1}$ of $X$
such that $Y$ isometrically embeds in $L^{p}(\Omega ,X_{1}).$ By \cite[Lemme
4.4]{Da1} $(iii)\Rightarrow (i)$ applied to $Z=L^{p}(\Omega ,X_{1}),$ it
suffices that $Z^{\ast }$ is G\^{a}teaux differentiable. Observe by \cite[%
Lemme 4.4]{Da1} ($(i)\Rightarrow (ii),$ that $(X_{1})^{\ast }$ is G\^{a}%
teaux differentiable). Thus we can assume that $X$ is separable. By \cite[%
p.~349]{Blas}, \cite[Chap~II-13-3, Corollary~1]{Din} $(L^{p}(\Omega
,X))^{\ast }=VB^{p^{\prime }}(\Omega ,X^{\ast }),$ where $p^{\prime }$ is
the conjugate of $p.$

\bigskip Let $f\in VB^{p^{\prime }}(\Omega ,X^{\ast })-\left\{ 0\right\} ;$
by lemma \ref{vvv}, it suffices to show that lim$_{R\ni t\rightarrow 0}\frac{%
\bigl\Vert f+th\bigr\Vert^{p^{\prime }}+\bigl\Vert f-th\bigr\Vert^{p^{\prime
}}-2\bigl\Vert f\bigr\Vert^{p^{\prime }}}{t}=0,$ for all $h\in VB^{p^{\prime
}}(\Omega ,X^{\ast }).$ Let $t>0$ and let $h\in VB^{p^{\prime }}(\Omega
,X^{\ast });$ we have

\begin{eqnarray*}
&&\text{lim}_{R\ni t\rightarrow 0}\frac{\bigl\Vert f+th\bigr\Vert%
_{VB^{p^{\prime }}(\Omega ,X^{\ast })}^{p^{\prime }}+\bigl\Vert f-th%
\bigr\Vert_{VB^{p^{\prime }}(\Omega ,X^{\ast })}^{p^{\prime }}-2\bigl\Vert f%
\bigr\Vert_{VB^{p^{\prime }}(\Omega ,X^{\ast })}^{p^{\prime }}}{t} \\
&=&\mathop{\displaystyle \int}\limits_{\Omega }\frac{\bigl\Vert f(\omega
)+th(\omega )\bigr\Vert_{X^{\ast }}^{p^{\prime }}+\bigl\Vert f(\omega
)-th(\omega )\bigr\Vert_{X^{\ast }}^{p^{\prime }}-2\bigl\Vert f(\omega )%
\bigr\Vert_{X^{\ast }}^{p^{\prime }}}{t}d\mu (\omega ) \\
&=&\mathop{\displaystyle \int}\limits_{\left\{ \omega ;\text{ }f(\omega
)=0\right\} }\frac{\bigl\Vert f(\omega )+th(\omega )\bigr\Vert_{X^{\ast
}}^{p^{\prime }}+\bigl\Vert f(\omega )-th(\omega )\bigr\Vert_{X^{\ast
}}^{p^{\prime }}-2\bigl\Vert f(\omega )\bigr\Vert_{X^{\ast }}^{p^{\prime }}}{%
t}d\mu (\omega )+ \\
&&\mathop{\displaystyle \int}\limits_{\left\{ \omega ;\text{ }f(\omega )\neq
0\right\} }\frac{\bigl\Vert f(\omega )+th(\omega )\bigr\Vert_{X^{\ast
}}^{p^{\prime }}+\bigl\Vert f(\omega )-th(\omega )\bigr\Vert_{X^{\ast
}}^{p^{\prime }}-2\bigl\Vert f(\omega )\bigr\Vert_{X^{\ast }}^{p^{\prime }}}{%
t}d\mu (\omega ) \\
&=&\mathop{\displaystyle \int}\limits_{\left\{ \omega ;\text{ }f(\omega
)=0\right\} }\frac{2t^{p^{\prime }}\bigl\Vert h(\omega )\bigr\Vert_{X^{\ast
}}^{p^{\prime }}}{t}d\mu (\omega )+ \\
&&\mathop{\displaystyle \int}\limits_{\left\{ \omega ;\text{ }f(\omega )\neq
0\right\} }\frac{\bigl\Vert f(\omega )+th(\omega )\bigr\Vert_{X^{\ast
}}^{p^{\prime }}+\bigl\Vert f(\omega )-th(\omega )\bigr\Vert_{X^{\ast
}}^{p^{\prime }}-2\bigl\Vert f(\omega )\bigr\Vert_{X^{\ast }}^{p^{\prime }}}{%
t}d\mu (\omega ).
\end{eqnarray*}

It is obvious that lim$_{R\ni t\rightarrow 0}$ $\mathop{\displaystyle \int}%
\limits_{\left\{ \omega ;\text{ }f(\omega )=0\right\} }\frac{2t^{p^{\prime }}%
\bigl\Vert h(\omega )\bigr\Vert_{X^{\ast }}^{p^{\prime }}}{t}d\mu (\omega
)=0.$ On the other hand by lemma \ref{vvv}, lim$_{R\ni t\rightarrow 0}\frac{%
\bigl\Vert f(\omega )+th(\omega )\bigr\Vert_{X^{\ast }}^{p^{\prime }}+%
\bigl\Vert f(\omega )-th(\omega )\bigr\Vert_{X^{\ast }}^{p^{\prime }}-2%
\bigl\Vert f(\omega )\bigr\Vert_{X^{\ast }}^{p^{\prime }}}{t}=0.$ By using
lemma \ref{UN} and the Lebesgue convergence theorem, we obtain%
\begin{equation*}
lim_{R\ni t\rightarrow 0}\mathop{\displaystyle \int}\limits_{\left\{ \omega ;%
\text{ }f(\omega )\neq 0\right\} }\frac{\bigl\Vert f(\omega )+th(\omega )%
\bigr\Vert_{X^{\ast }}^{p^{\prime }}+\bigl\Vert f(\omega )-th(\omega )%
\bigr\Vert_{X^{\ast }}^{p^{\prime }}-2\bigl\Vert f(\omega )\bigr\Vert%
_{X^{\ast }}^{p^{\prime }}}{t}d\mu (\omega )=0.
\end{equation*}

This implies that lim$_{R\ni t\rightarrow 0}\mathop{\displaystyle \int}%
\limits_{\Omega }\frac{\bigl\Vert f(\omega )+th(\omega )\bigr\Vert_{X^{\ast
}}^{^{\prime }p}+\bigl\Vert f(\omega )-th(\omega )\bigr\Vert_{X^{\ast
}}^{p^{\prime }}-2\bigl\Vert f(\omega )\bigr\Vert_{X^{\ast }}^{p^{\prime }}}{%
t}d\mu (\omega )=0.$ \qed

\end{document}